\newcommand{\CLASSINPUTtoptextmargin}{54pt}
\newcommand{\CLASSINPUTbottomtextmargin}{72pt}
\newcommand{\CLASSINPUTinnersidemargin}{45pt}
\newcommand{\CLASSINPUToutersidemargin}{45pt}
\documentclass[a4paper, conference]{IEEEtran}
\IEEEoverridecommandlockouts
\usepackage{cite}
\usepackage{amsmath,amssymb,amsfonts,amsthm}
\usepackage{algorithmic}
\usepackage{graphicx}
\usepackage{textcomp}
\usepackage{xcolor}
\usepackage{enumitem}
\usepackage{url}

\DeclareMathOperator{\Real}{Re}

\theoremstyle{plain}

\theoremstyle{remark}
\newtheorem{remark}{Remark}

\begin{document}

\urlstyle{tt}

\setlist[enumerate]{label={(\alph*)}, ref={(\alph*)}, leftmargin=*}
\newlist{steps}{enumerate}{1}
\setlist[steps]{label={\arabic*.}, ref={\arabic*}, leftmargin=*}

\title{Partial Pole Placement via Delay Action: A Python Software for Delayed Feedback Stabilizing Design\thanks{This work is partially supported by a public grant overseen by the French National Research Agency (ANR) as part of the ``Investissement d'Avenir'' program, through the iCODE project funded by the IDEX Paris-Saclay, ANR-11-IDEX0003-02. The authors also acknowledge the support of Institut Polytechnique des Sciences Avanc\'ees (IPSA).}}

\author{\IEEEauthorblockN{Islam Boussaada\IEEEauthorrefmark{1}\IEEEauthorrefmark{2}, Guilherme Mazanti\IEEEauthorrefmark{1}\IEEEauthorrefmark{2}, Silviu-Iulian Niculescu\IEEEauthorrefmark{1},\\Julien Huynh\IEEEauthorrefmark{2}\IEEEauthorrefmark{3}, Franck Sim\IEEEauthorrefmark{2}\IEEEauthorrefmark{3}, Matthieu Thomas\IEEEauthorrefmark{2}\IEEEauthorrefmark{3}}\IEEEauthorblockA{\IEEEauthorrefmark{1}Universit\'e Paris-Saclay, CNRS, CentraleSup\'elec, Inria \\ Laboratoire des signaux et syst\`emes (L2S) \\ 91190 Gif-sur-Yvette, France \\ E-mails: \{first name.last name\}@l2s.centralesupelec.fr}\IEEEauthorblockA{\IEEEauthorrefmark{2}Institut Polytechnique des Sciences Avanc\'ees (IPSA) \\ 63 boulevard de Brandebourg, 94200 Ivry-sur-Seine, France \\ E-mails: \{first name.last name\}@ipsa.fr}\IEEEauthorblockA{\IEEEauthorrefmark{3}Cyb'Air Association, 94200 Ivry-sur-Seine, France}}

\maketitle

\begin{abstract}
This paper presents a new Python software for the parametric design of stabilizing feedback laws with time delays, called \emph{Partial Pole Placement via Delay Action} (P3$\pmb\delta$). After an introduction recalling recent theoretical results on the multiplicity-induced-dominancy (MID) and coexisting real roots-induced-dominancy (CRRID) properties and their use for the feedback stabilization of control systems operating under time delays, the paper presents the current version of P3$\pmb\delta$, which relies on the MID property to compute delayed stabilizing feedback laws for scalar differential equations with a single delay. We detail in particular its graphical user interface (GUI), which allows the user to input the necessary information and obtain the results of the analysis done by the software. These results include the parameters stabilizing the closed-loop system, graphical representations of the spectrum of the closed-loop system, simulations of solutions in the time domain, and a sensitivity analysis with respect to uncertain delays.
\end{abstract}

\begin{IEEEkeywords}
Time-delay systems, Controller design, Stability, Stabilization, Python toolbox, GUI
\end{IEEEkeywords}

\section{Introduction}

Time delays often occur in control systems, mainly due to the time required for acquiring, propagating, or processing information. For this reason, systems with time delays are a frequent topic in the control theory literature, with many works, such as \cite{Gu2003Stability, Hale1993Introduction, Michiels2014Stability}, highlighting the effects of delays on the behavior of control systems, in particular on their stability.

Commonly, time delays lead to desynchronizing or destabilizing effects on the dynamics of the system they appear. However, some works have emphasized that the delay may also have a stabilizing effect in control design. For instance, in \cite{Tallman1958Analog}, a delayed controller is used in order to improve the stability of systems with oscillatory behavior and small damping. The stabilization properties of delayed controllers has also been considered in \cite{Suh1979Proportional}, which uses a proportional-delayed controller, replacing the classical proportional-derivative controller thanks to the ``average derivative action'' obtained via the time delay, a technique also used in \cite{Atay1999Balancing}. Further discussion of the stabilizing effects of time delays can be found in \cite{Niculescu2010Delay}, which highlights in particular the fact that closed-loop stability may be guaranteed for some control systems precisely by the existence of the delay. A growing literature exhibits the design of delayed controllers in a wide range of applications, such as, for instance, the control of flexible mechanical structures or the regulation of networks (see, e.g., \cite{Boussaada2018Further, Irofti2016Codimension}).

In this paper, we consider linear time-invariant differential equations with a single time delay under the form
\begin{multline}
\label{MainEqn}
y^{(n)}(t) + a_{n-1} y^{(n-1)}(t) + \dotsb + a_0 y(t) \\ + b_m y^{(m)}(t - \tau) + \dotsb + b_0 y(t - \tau) = 0,
\end{multline}
where $\tau > 0$ is the positive delay, $y$ is the real-valued unknown function, $n$ and $m$ are nonnegative integers with $n > m$, and $a_0, \dotsc, a_{n-1}, b_0, \dotsc, b_m$ are real coefficients.

The stability analysis of a linear time-invariant time-delay system can be addressed using spectral methods by considering the corresponding characteristic function, whose complex roots determine the asymptotic behavior of solutions of the system, as presented, e.g., in \cite{Hale1993Introduction, Michiels2014Stability}. The characteristic function corresponding to \eqref{MainEqn} is
\begin{equation}
\label{Delta}
\Delta(s) = s^n + \sum_{k=0}^{n-1} a_k s^k + e^{-s \tau} \sum_{k=0}^m b_k s^k,
\end{equation}
and \eqref{MainEqn} is exponentially stable if and only if the \emph{spectral abscissa} $\gamma = \sup\{\Real s \mid \Delta(s) = 0\}$ satisfies $\gamma < 0$. Equation \eqref{MainEqn} is said to be of \emph{retarded type}, since the highest-order derivative only appears in the non-delayed term $y^{(n)}(t)$.

Equations under the form \eqref{MainEqn} may arise from linear time-invariant controlled differential equations, such as $y^{(n)} + a_{n-1} y^{(n-1)}(t) + \dotsb + a_0 y(t) = u(t)$, when applying a delayed feedback control under the form $u(t) = - b_m y^{(m)}(t - \tau) - \dotsb - b_0 y(t - \tau)$. In this case, the behavior of the closed-loop system is the influenced by the choices of the free parameters $b_0, \dotsc, b_m$ in the feedback control, which are thus free coefficients in the characteristic function \eqref{Delta}.

The characteristic function \eqref{Delta} is a particular case of a \emph{quasipolynomial}, i.e., a polynomial in the variables $s$ and $e^{-s \tau}$. Quasipolynomials have been considered in several works, such as \cite{Berenstein1995Complex, Hale1993Introduction, Stepan1989Retarded, Wielonsky2001Rolle}, often in connection with the analysis of time-delay systems. A major difficulty in the study of quasipolynomials for the feedback stabilization of time-delay systems is that quasipolynomials have infinitely many roots, but one only disposes of finitely many parameters in the feedback law to choose the location of these roots and place them in order to guarantee a negative spectral abscissa, and hence exponential stability of the closed-loop system.

\newcommand{\DominancyRefs}{\cite{Amrane2018Qualitative, BedouheneReal, Boussaada2018Dominancy, Boussaada2020Multiplicity, Boussaada2018Further, MazantiMultiplicity, Mazanti2020Qualitative, Mazanti2020Spectral}}

Recent works such as \DominancyRefs{} have been interested in the design of pole placement techniques for quasipolynomials with the aim of selecting the free parameters of the system in order to choose the location of finitely many roots in the complex plane and guarantee that the \emph{dominant root}, i.e., the rightmost root on the complex plane, is among the chosen ones. Unlike methods based on finite spectrum assignment such as those from \cite{Manitius1979Finite}, the controllers designed using these techniques do not render the closed-loop system finite dimensional, but control instead its rightmost spectral value. 

The works \DominancyRefs{} usually proceed either by assigning a real root of maximal multiplicity and proving that this root is necessarily the rightmost root of the characteristic quasipolynomial (a property known as \emph{multiplicity-induced-dominancy}, or MID for short) or by assigning a certain amount of real roots (typically equally spaced for simplicity) and proving that the rightmost root among the assigned roots is also the rightmost root of the characteristic quasipolynomial (a property known as \emph{coexisting real roots-induced-dominancy}, or CRRID for short).

The MID property for \eqref{MainEqn} is shown, for instance, in \cite{Boussaada2018Further} in the case $n = 2$ and $m = 0$, in \cite{Boussaada2020Multiplicity} in the case $n = 2$ and $m = 1$ (see also \cite{Boussaada2018Dominancy}), and in \cite{MazantiMultiplicity} in the case of any positive integer $n$ and $m = n-1$ (see also \cite{Mazanti2020Qualitative}). The CRRID property is shown, for instance, in \cite{Amrane2018Qualitative} in the cases $(n, m) = (2, 0)$ and $(n, m) = (1, 0)$, and in \cite{BedouheneReal} in the case of any positive integer $n$ and $m = 0$. In all these cases, the maximal multiplicity of a real root or, equivalently, the maximal number of coexisting simple real roots is the integer $n + m + 1$.

This paper presents the \emph{Partial Pole Placement via Delay Action} software (P3$\delta$ for short), a Python software based on the results from \DominancyRefs{} for the parametric design of stabilizing feedback laws with time delays. The first version of P3$\delta$, presented in the current paper, allows for the design of feedback laws for linear time-invariant differential equations with a single time delay under the form \eqref{MainEqn} using MID techniques.

Several other softwares have been recently developed for the analysis of time-delay systems from various perspectives, such as stability, robustness, or bifurcation aspects. This is the case, for instance, of the Matlab packages YALTA \cite{Avanessoff2014Hinfty}, dedicated to the $H_\infty$ stability analysis of time-delay systems with commensurate delays, TRACE-DDE \cite{Breda2009Trace}, devoted to the computation of characteristic roots and stability charts of linear autonomous time-delay systems, DDE-BIFTOOL \cite{Engelborghs2002Numerical}, interested in the computation, continuation, and stability analysis of steady-state solutions of time-delay systems and their bifurcations, and QPmR \cite{Vyhlidal2014QPmR}, specialized in the computation of roots of quasipolynomials. One of the major novelties of P3$\delta$ lies in addressing the stabilization of control systems with time-delays by using of the MID property to design stabilizing feedback laws. For that purpose, P3$\delta$ makes use of both symbolic and numeric computations.

\section{Description of P3$\delta$}
\label{SecP3Delta}

P3$\delta$ is freely available for download on \url{https://cutt.ly/p3delta}, where installation instructions, video demonstrations, and the user guide are also available. Interested readers may also contact directly any of the authors of the paper.



In the current version of P3$\delta$, only the MID property is exploited for the stabilization of \eqref{MainEqn}. This can be done in two different ways, named ``Classic MID'' and ``Control-oriented MID'', according to which coefficients of \eqref{MainEqn} are assumed to be fixed and which are assumed to be free.


\subsection{Classic MID mode}
\label{SecClassicMID}

The ``Classic MID'' mode corresponds to considering that all coefficients of the quasipolynomial $\Delta$ from \eqref{Delta} are free. The user inputs the values of the delay $\tau$ and of the desired real root $s_0$ and P3$\delta$ computes all coefficients $a_0, \dotsc, a_{n-1}, b_0, \dotsc, b_{m}$ ensuring that the value $s_0$ is a dominant root of $\Delta$ of maximal multiplicity $n + m + 1$. To use the ``Classic MID'' mode, the user should proceed as follows:
\begin{steps}
\item Enter the values of the integers $n$ and $m$ appearing in the differential equation \eqref{MainEqn}.
\item Select the ``Classic MID'' option in the drop-down menu ``--- Choose MID type ---''.
\end{steps}
After this selection, the window of the program is filled with the places for the other inputs and the outputs of P3$\delta$.
\begin{steps}[resume]
\item Enter the values of the desired real root of maximal multiplicity $s_0$ and of the delay $\tau$ in the corresponding fields that appear below the drop-down menu.
\item\label{ClassicMIDStepConfirm} Enter the bounds $x_{\min}, x_{\max}, y_{\min}, y_{\max}$ of the rectangle $[x_{\min}, x_{\max}] \times [y_{\min}, y_{\max}] \subset \mathbb C$ in which P3$\delta$ will look for roots of \eqref{Delta} and press the ``Confirm'' button.
\end{steps}
Once the ``Confirm'' button is pressed, P3$\delta$ will compute the values of the coefficients $a_0, \dotsc, a_{n-1}, b_0, \dotsc, b_m$ ensuring that $s_0$ is a root of maximal multiplicity of the quasipolynomial $\Delta$ from \eqref{Delta} and show their values. P3$\delta$ will also numerically compute all roots of $\Delta$ within the selected rectangle using the computed values of $a_0, \dotsc, a_{n-1},\allowbreak b_0, \dotsc, b_m$ and plot these roots in the plot ``Roots'' at the lower left corner of the window. This numerical computation is carried out using Python's \texttt{cxroots} module, which implements numerical methods described in \cite{Kravanja2000Computing}. 

Optionally, after the previous computations are completed, the user may also simulate some trajectories of the system in the time domain. This can be done, after completing step \ref{ClassicMIDStepConfirm} above, by the following steps:
\begin{steps}[resume]
\item\label{ClassicMIDStepTimeSimulationBegin} Choose the type of the initial condition from the drop-down menu ``--- Initial Solution ---''.
\end{steps}
The currently supported types are ``Constant'', ``Polynomial'', ``Exponential'', and ``Trigonometric'', which corresponds to initial conditions of the forms $x(t) = c$, $x(t) = \sum_{k=0}^r c_k t^k$, $x(t) = A e^{\gamma t}$, and $x(t) = A \sin(\omega t + \varphi)$, respectively, where $c, r, c_0, \dotsc, c_r, A, \gamma, \omega, \varphi$ are constants to be chosen by the user and the initial condition is defined in the time interval $[-\tau, 0]$.
\begin{steps}[resume]
\item\label{ClassicMIDStepT} Enter the simulation time in the corresponding box.
\item Enter the values of the constants appearing in the expression of the initial condition in the corresponding input boxes.
\item\label{ClassicMIDStepTimeSimulationEnd} After entering all the constants, press ``Enter'' on the keyboard or click on the ``Confirm'' button appearing in the same frame as the constants.
\end{steps}
After these steps, the numerical solution corresponding to the chosen initial condition will be computed using an explicit Euler scheme in the time interval $[-\tau, T]$, where $T$ is the value entered in step \ref{ClassicMIDStepT}. The corresponding solution will be plotted in the graph on the ``Solutions'' part of the screen.

\begin{figure*}[ht]
\centering
\includegraphics[width=0.8\textwidth]{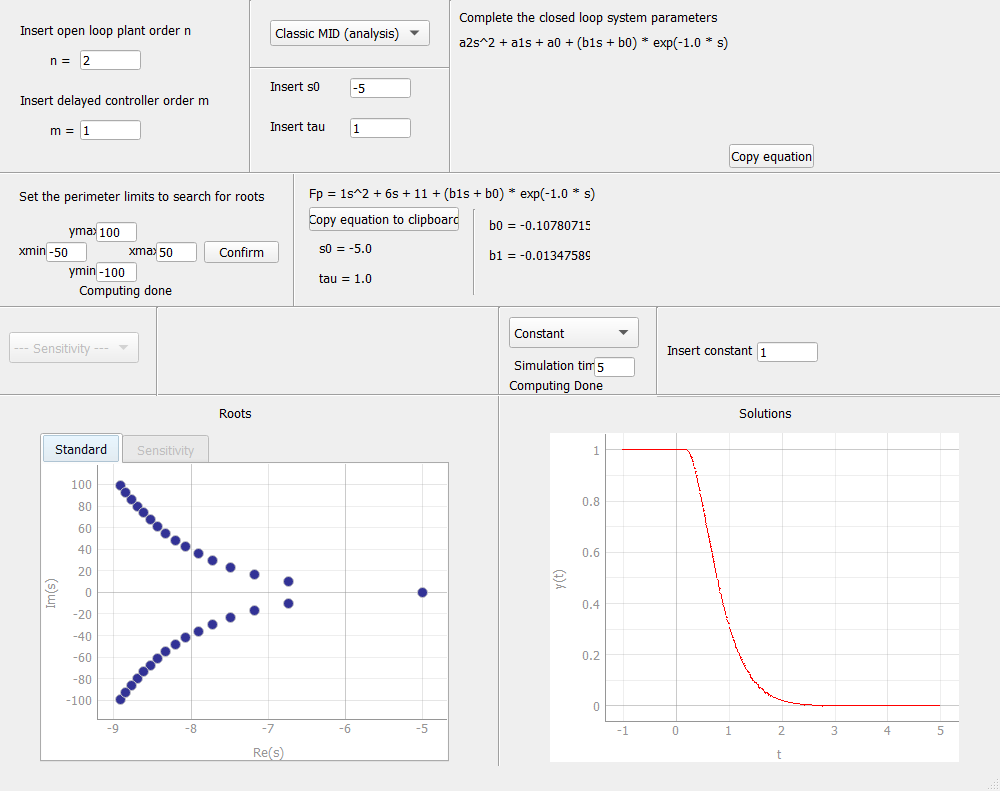}
\caption{``Classic MID'' mode of P3$\delta$.}
\label{FigClassicMID}
\end{figure*}

Figure~\ref{FigClassicMID} shows a screen capture of the ``Classic MID'' mode of P3$\delta$. In this figure, we have chosen $n = 2$, $m = 1$, the ``Classic MID'' mode, and the values $s_0 = -5$ and $\tau = 1$. After entering $x_{\min} = -50$, $x_{\max} = 50$, $y_{\min} = -100$, and $y_{\max} = 100$ and having clicked on ``Confirm'', P3$\delta$ shows the values of the coefficients of the quasipolynomial, $a_1 = 6$, $a_0 = 11$, $b_1 \approx -0.1078$, and $b_0 \approx -0.0135$, ensuring that $s_0 = -5$ is a root of maximal multiplicity $n + m + 1 = 4$ in this case. P3$\delta$ also plots the numerical roots of the quasipolynomial in the selected rectangle in the graph on the lower left corner. After choosing the ``Constant'' initial condition, selecting the simulation time $T = 5$ and the value $1$ for the constant, the corresponding numerical solution of the system appears in the graph in the lower right corner of the screen.

\subsection{Control-oriented MID mode}
\label{SecControlMID}

The ``Control-oriented MID'' mode corresponds to considering that the coefficients $a_0, \dotsc, a_{n-1}$ corresponding to the non-delayed terms of \eqref{MainEqn} are given and that the coefficients $b_0, \dotsc, b_m$ corresponding to the delayed terms are free. The user may choose to input either the value of $\tau$ or that of $s_0$ (but not both) and P3$\delta$ computes all coefficients $b_0, \dotsc, b_m$ ensuring the existence of a dominant root of the quasipolynomial $\Delta$ from \eqref{Delta} of multiplicity $m+2$. P3$\delta$ also computes the value of the parameter among $\tau$ or $s_0$ that has not been fixed by the user.

\begin{remark}
In the ``Control-oriented MID'' mode, it may happen to be impossible to choose a real root $s_0$ of multiplicity $m+2$. In this case, P3$\delta$ warns the user of this fact and provides an equation relating $s_0$ and $\tau$. The user should either enter a value of $s_0$ such that this equation admits a positive root $\tau$ or a positive value of $\tau$ such that this equation admits a real root $s_0$ in order to proceed with the computations.
\end{remark}

To use the ``Control-oriented MID'' mode, the user should proceed as follows:
\begin{steps}
\item Enter the values of the integers $n$ and $m$ appearing in the differential equation \eqref{MainEqn}.
\item Select the ``Control-oriented MID'' option in the drop-down menu ``--- Choose MID type ---''.
\end{steps}
After this selection, the window of the program is filled with the places for the other inputs and the outputs of P3$\delta$.
\begin{steps}[resume]
\item\label{ControlOrientedMIDStepChooses0Tau} Select from the drop-down menu ``--- Choose s0 or tau ---'' whether to input the value of the multiple root $s_0$ or the value of the delay $\tau$.
\item Enter the value of $s_0$ or $\tau$, according to the choice of the previous step.
\item Enter the values of the known coefficients $a_0, \dotsc, a_{n-1}$ and press the ``Confirm'' button located in the same frame.
\item Enter the bounds $x_{\min}, x_{\max}, y_{\min}, y_{\max}$ of the rectangle $[x_{\min}, x_{\max}] \times [y_{\min}, y_{\max}] \subset \mathbb C$ in which P3$\delta$ will look for roots of \eqref{Delta}.
\item\label{ControlOrientedMIDStepConfirm} Press the ``Confirm'' button.
\end{steps}
Once the ``Confirm'' button is pressed, P3$\delta$ will compute the values of the coefficients $b_0, \dotsc, b_m$ ensuring that $s_0$ is a root of multiplicity $m+2$ of the quasipolynomial $\Delta$ from \eqref{Delta} and show their values. Similarly to the ``Classic MID'' option, P3$\delta$ will also numerically compute all roots of $\Delta$ within the selected rectangle by using Python's \texttt{cxroots} module and output the result in the plot ``Roots''.

As in the ``Classic MID'' case, the user may plot solutions in the time domain. After completing step \ref{ControlOrientedMIDStepConfirm} above, this can be done by following the same steps \ref{ClassicMIDStepTimeSimulationBegin}--\ref{ClassicMIDStepTimeSimulationEnd} from Section~\ref{SecClassicMID}.

\begin{figure*}[ht]
\centering
\includegraphics[width=0.8\textwidth]{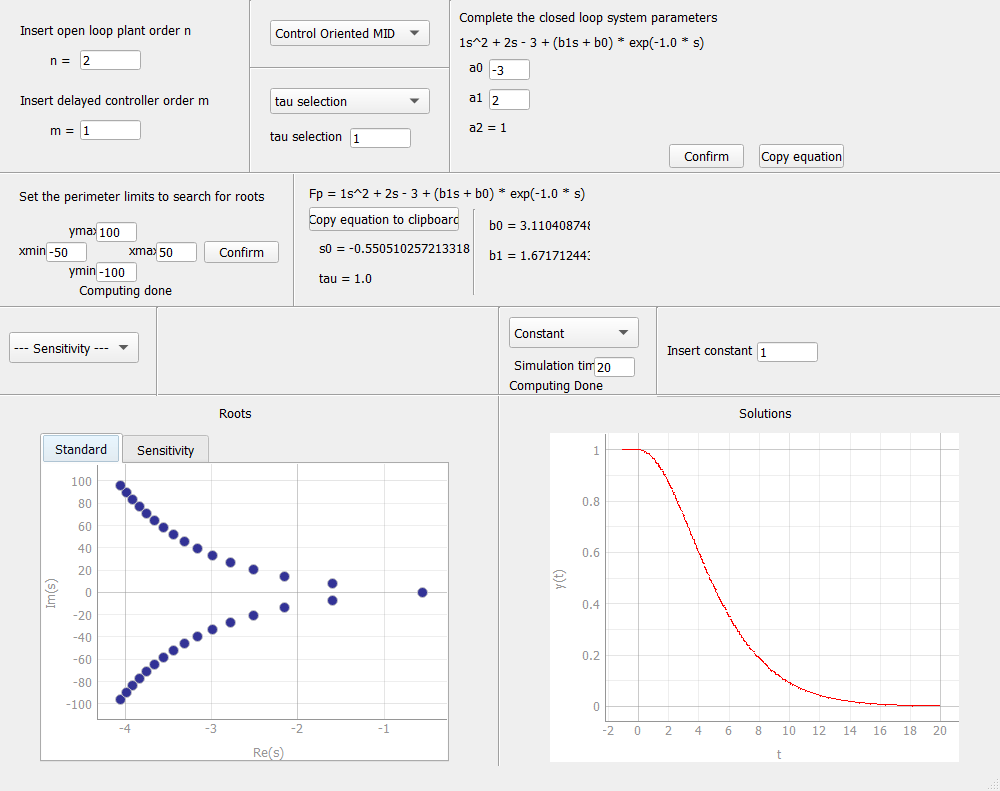}
\caption{``Control-oriented MID'' mode of P3$\delta$.}
\label{FigControlOrientedMID}
\end{figure*}

Figure~\ref{FigControlOrientedMID} shows a screen capture of the ``Control-oriented MID'' mode of P3$\delta$. In this figure, we have chosen $n = 2$, $m = 1$, the ``Control-oriented MID'' mode, the input of $\tau$, the value $\tau = 1$, and the coefficients $a_0 = -3$ and $a_1 = 2$. After entering $x_{\min} = -50$, $x_{\max} = 50$, $y_{\min} = -100$, and $y_{\max} = 100$ and having clicked on ``Confirm'', P3$\delta$ shows the values of the coefficients of the quasipolynomial, $b_1 \approx 1.6717$ and $b_0 \approx 3.1104$, recalls the value $\tau = 1.0$, and also shows the value of the root of multiplicity $m+2 = 3$, $s_0 \approx -0.550510$. P3$\delta$ also plots the numerical roots of the quasipolynomial in the selected rectangle in the graph on the lower left corner. After choosing the ``Constant'' initial condition, selecting the simulation time $T = 20$ and the value $1$ for the constant, the corresponding numerical solution of the system appears in the graph in the lower right corner of the screen.

In addition to these outputs, which are similar to the ``Classic MID'' case, the ``Control-oriented MID'' option can also perform a numerical sensitivity analysis of the computed roots with respect to variations in the delay $\tau$. To do so, the user should follow the above steps up to step \ref{ControlOrientedMIDStepConfirm}, selecting to enter the value of $\tau$ in step \ref{ControlOrientedMIDStepChooses0Tau}. Then, the steps to get the sensitivity plot are the following:
\begin{steps}
\item Select the ``Sensitivity'' tab in the ``Roots'' plot.
\item Select ``tau sensitivity'' in the drop-down menu ``--- Sensitivity ---'' above the ``Roots'' plot.
\item Enter the value of the step $\varepsilon$ and the number of iterations $K$ in the corresponding boxes.
\item Enter the bounds $x_{\min}, x_{\max}, y_{\min}, y_{\max}$ of the rectangle $[x_{\min}, x_{\max}] \times [y_{\min}, y_{\max}] \subset \mathbb C$ in which P3$\delta$ will look for roots of \eqref{Delta}.
\end{steps}
Since the sensitivity computation may take quite some time, it is highly recommended to choose a smaller rectangle containing few roots of $\Delta$, including the dominant multiple root.
\begin{steps}[resume]
\item Press the ``Confirm'' button in the frame of the bounds of the rectangle.
\end{steps}
Once these steps are completed, the sensitivity plot appears in the ``Roots'' plot. This plot contains the roots of $\Delta$ in the selected rectangle for the values of delays $\tau + k \varepsilon$ for $k \in \{-K, -K+1, \dotsc, K-1, K\}$. Roots computed with negative values of $k$, corresponding to values of the delay smaller than $\tau$, are represented in shades of blue, with darker blue representing $k = -K$ and lighter tones representing increasing values of $k$. Roots computed with positive values of $k$, corresponding to values of the delay larger than $\tau$, are represented in shades of orange to red, with darker red representing $k = K$ and lighter tones moving to orange representing decreasing values of $k$. The roots computed with $k = 0$, corresponding to the nominal value of $\tau$ selected by the user, are represented by black diamonds.

\begin{figure}[ht]
\centering
\includegraphics[width=0.8\columnwidth]{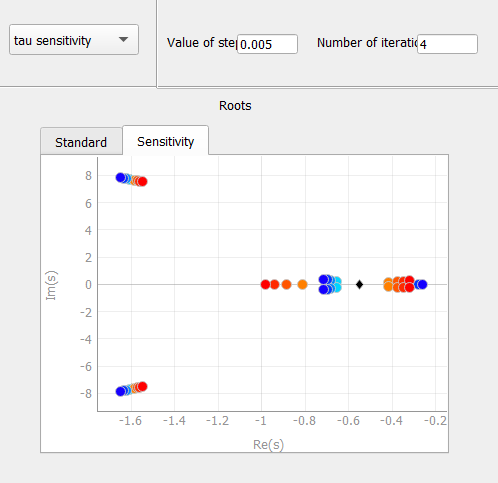}
\caption{Detail of the P3$\delta$ screen for computing sensitivity with respect to the delay.}
\label{FigSensitivity}
\end{figure}

Figure~\ref{FigSensitivity} represents the part of P3$\delta$ screen corresponding to the sensitivity computation. After having completed the steps that led to the screen shown in Figure~\ref{FigControlOrientedMID}, we have chosen the step $\varepsilon = 0.005$ and the number of steps $K = 4$ and selected the default rectangle $[-5, 5] \times [-10, 10] \subset \mathbb C$. Clicking on ``Confirm'', P3$\delta$ outputs the graph shown in Figure~\ref{FigSensitivity}. We observe that, when the delay $\tau$ is perturbed, the root $s_0 \approx -0.550510$ of multiplicity $3$ splits into three simple roots.

\section{Illustrative examples}
\label{SecExpl}

As illustrations of the use of P3$\delta$, this section revisits two examples from \cite{Boussaada2020Multiplicity}.

\subsection{A first order equation}

We consider here the delay-differential equation
\begin{equation}
\label{EqExpl1}
\dot y(t) + a_0 y(t) + b_0 y(t - \tau) = 0,
\end{equation}
whose characteristic quasipolynomial is $\Delta(s) = s + a_0 + b_0 e^{-s \tau}$. We then have $n = 1$ and $m = 0$. According to \cite{Boussaada2020Multiplicity}, the maximal multiplicity of a real root of $\Delta$ is $2$, and it is attained if and only if $a_0 = -s_0 - \frac{1}{\tau}$ and $b_0 = \frac{e^{s_0 \tau}}{\tau}$.

Inputting $n = 1$, $m = 0$, selecting ``Classic MID'', and choosing $s_0 = -2$ and $\tau = 1$ in P3$\delta$, we obtain $a_0 = 1$ and $b_0 \approx 0.1353$, which is in accordance with the above expressions for $a_0$ and $b_0$. We also obtain the roots of $\Delta$ in a given rectangle, represented in Figure~\ref{FigExpl}(a) for the rectangle $[-50, 50] \times [-100, 100]$, and time simulations of solutions, for instance the one from Figure~\ref{FigExpl}(b), obtained with initial condition $x(t) = \sin(10 t)$.

\subsection{Stabilization of the double integrator}

Let us consider a double integrator $\ddot y(t) = u(t)$ with the delayed feedback control $u(t) = - b_1 \dot y(t - \tau) - b_0 y(t - \tau)$, which yields the delay-differential equation
\begin{equation}
\label{EqExpl2}
\ddot y(t) + b_1 \dot y(t-\tau) + b_0 y(t - \tau) = 0,
\end{equation}
whose characteristic quasipolynomial is $\Delta(s) = s^2 + (b_1 s + b_0) e^{-s \tau}$. This corresponds to $n = 2$, $m = 1$, and $a_0 = a_1 = 0$. According to \cite{Boussaada2020Multiplicity}, the maximal achievable multiplicity for a root $s_0$ of $\Delta$ is $3$, which is attained if and only if
\begin{equation}
\label{Expl2Conditions}
\begin{gathered}
b_1 = 2 \frac{\left(\sqrt{2} - 1\right) e^{- 2 + \sqrt{2}}}{\tau}, \quad b_0 = 2 \frac{\left(5\sqrt{2} - 7\right) e^{- 2 + \sqrt{2}}}{\tau^2}, \\
s_0 = -\frac{2 - \sqrt{2}}{\tau}.
\end{gathered}
\end{equation}

\begin{figure*}
\centering
\begin{tabular}{cc}
\includegraphics[width=0.4\textwidth]{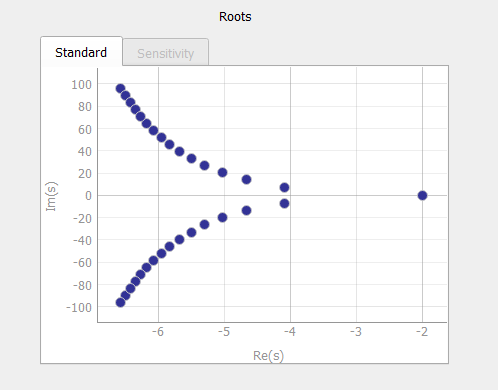} & \includegraphics[width=0.4\textwidth]{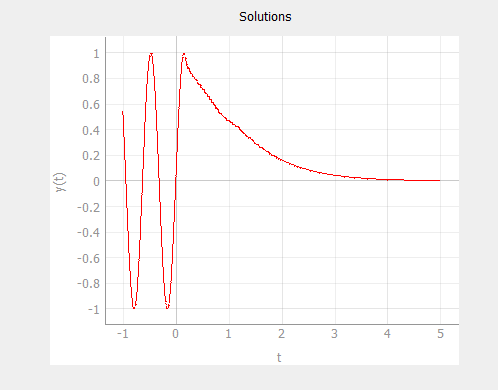} \\
(a) & (b) \\
\includegraphics[width=0.4\textwidth]{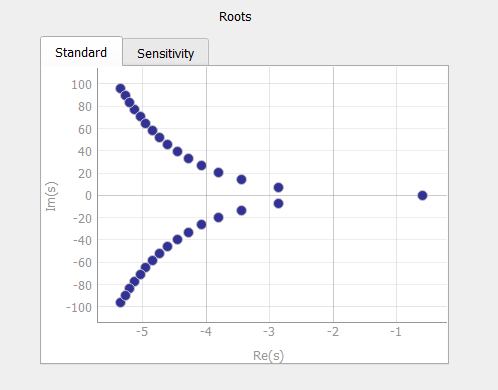} & \includegraphics[width=0.4\textwidth]{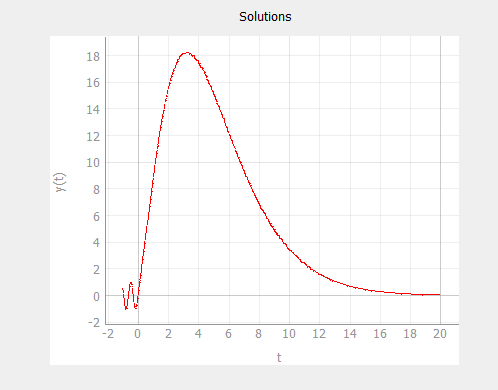} \\
(c) & (d)
\end{tabular}
\caption{Example of the usage of P3$\delta$ (a, b) in the ``Classic MID'' mode for the delay-differential equation \eqref{EqExpl1} and (c, d) in the ``Control-oriented MID'' mode for the delay-differential equation \eqref{EqExpl2}. (a, c) Roots of the characteristic quasipolynomial in the rectangle $[-50, 50] \times [-100, 100]$. (b, d) Numerical simulation of a solution with initial condition $x(t) = \sin(10 t)$.}
\label{FigExpl}
\end{figure*}

Inputting $n = 2$, $m = 1$, selecting ``Control-oriented MID'', and choosing $\tau = 1$, $a_0 = 0$, and $a_1 = 0$ in P3$\delta$, we obtain $b_1 \approx 0.4612$, $b_0 \approx 0.0791$, and $s_0 \approx -0.585786$, which is in accordance with \eqref{Expl2Conditions}. We also obtain the roots of $\Delta$ in a given rectangle, represented in Figure~\ref{FigExpl}(c) for the rectangle $[-50, 50] \times [-100, 100]$, and time simulations of solutions, for instance the one from Figure~\ref{FigExpl}(d), obtained with initial condition $x(t) = \sin(10 t)$.

\section{Conclusion and planned developments}

By using the recent theoretical results on the MID property for linear delay-differential equations under the form \eqref{MainEqn}, P3$\delta$ computes the system parameters ensuring that a given real root $s_0$ attains its maximal multiplicity and is thus dominant. P3$\delta$ currently works in two modes, ``Classic MID'' and ``Control-oriented MID'', which differ on which coefficients of \eqref{MainEqn} are assumed to be free. In all modes, P3$\delta$ computes the values of the free coefficients ensuring maximal multiplicity, performs a numerical computation of the roots of the characteristic quasipolynomial, and is able to perform time-domain simulations. The ``Control-oriented MID'' option also offers the possibility of a numerical study of the sensitivity of the roots with respect to variations in the delay.

There are currently several plans for future development of P3$\delta$ in its future versions, including sensitivity analysis with respect to the known parameters $a_0, \dotsc, a_{n-1}$ in the ``Control-oriented MID'' option, the inclusion of new options allowing for a more flexible choice of which coefficients are assumed to be fixed and which are assumed to be free, and the use of CRRID-based results to place coexisting real roots.

\section*{Acknowledgments}

The authors wish to acknowledge the work of the full P3$\delta$ development team, which, in addition to the authors, also include Mickael Alcaniz, Yoann Audet, Thomas Charbonnet, Honor\'{e} Curlier, Adrien Leclerc, Max Perraudin, Pierre-Henry Poret, and Achrafy Said Mohamed. The development of P3$\delta$ was also made possible thanks to the work of the Cyb'Air Association.

\bibliographystyle{IEEEtran}
\bibliography{p3delta}

\begin{thebibliography}{10}
\providecommand{\url}[1]{#1}
\csname url@samestyle\endcsname
\providecommand{\newblock}{\relax}
\providecommand{\bibinfo}[2]{#2}
\providecommand{\BIBentrySTDinterwordspacing}{\spaceskip=0pt\relax}
\providecommand{\BIBentryALTinterwordstretchfactor}{4}
\providecommand{\BIBentryALTinterwordspacing}{\spaceskip=\fontdimen2\font plus
\BIBentryALTinterwordstretchfactor\fontdimen3\font minus
  \fontdimen4\font\relax}
\providecommand{\BIBforeignlanguage}[2]{{%
\expandafter\ifx\csname l@#1\endcsname\relax
\typeout{** WARNING: IEEEtran.bst: No hyphenation pattern has been}%
\typeout{** loaded for the language `#1'. Using the pattern for}%
\typeout{** the default language instead.}%
\else
\language=\csname l@#1\endcsname
\fi
#2}}
\providecommand{\BIBdecl}{\relax}
\BIBdecl

\bibitem{Gu2003Stability}
K.~Gu, V.~L. Kharitonov, and J.~Chen, \emph{Stability of time-delay systems},
  ser. Control Engineering.\hskip 1em plus 0.5em minus 0.4em\relax
  Birkh\"{a}user Boston, Inc., Boston, MA, 2003.

\bibitem{Hale1993Introduction}
J.~K. Hale and S.~M. Verduyn~Lunel,
  \emph{\BIBforeignlanguage{english}{Introduction to functional differential
  equations}}.\hskip 1em plus 0.5em minus 0.4em\relax New York:
  Springer-Verlag, 1993, vol.~99.

\bibitem{Michiels2014Stability}
W.~Michiels and S.-I. Niculescu, \emph{Stability, control, and computation for
  time-delay systems: An eigenvalue-based approach}, 2nd~ed.\hskip 1em plus
  0.5em minus 0.4em\relax SIAM, Philadelphia, PA, 2014.

\bibitem{Tallman1958Analog}
G.~H. Tallman and O.~J.~M. Smith, ``Analog study of dead-beat posicast
  control,'' \emph{IRE Transactions on Automatic Control}, vol.~4, no.~1, pp.
  14--21, 1958.

\bibitem{Suh1979Proportional}
I.~H. Suh and Z.~Bien, ``Proportional minus delay controller,'' \emph{IEEE
  Trans. Automat. Control}, vol.~24, no.~2, pp. 370--372, 1979.

\bibitem{Atay1999Balancing}
F.~M. Atay, ``Balancing the inverted pendulum using position feedback,''
  \emph{Appl. Math. Lett.}, vol.~12, no.~5, pp. 51--56, 1999.

\bibitem{Niculescu2010Delay}
S.-I. Niculescu, W.~Michiels, K.~Gu, and C.~T. Abdallah, ``Delay effects on
  output feedback control of dynamical systems,'' in \emph{Complex time-delay
  systems}, F.~M. Atay, Ed.\hskip 1em plus 0.5em minus 0.4em\relax Springer,
  Berlin, 2010, pp. 63--84.

\bibitem{Boussaada2018Further}
I.~Boussaada, S.~Tliba, S.-I. Niculescu, H.~U. \"{U}nal, and T.~Vyhl\'{\i}dal,
  ``Further remarks on the effect of multiple spectral values on the dynamics
  of time-delay systems. {A}pplication to the control of a mechanical system,''
  \emph{Linear Algebra Appl.}, vol. 542, pp. 589--604, 2018.

\bibitem{Irofti2016Codimension}
D.-A. Irofti, I.~Boussaada, and S.-I. Niculescu, ``On the codimension of the
  singularity at the origin for networked delay systems,'' in \emph{Delays and
  Networked Control Systems}, A.~Seuret, L.~Hetel, J.~Daafouz, and K.~H.
  Johansson, Eds.\hskip 1em plus 0.5em minus 0.4em\relax Springer International
  Publishing, 2016, pp. 3--15.

\bibitem{Berenstein1995Complex}
C.~A. Berenstein and R.~Gay, \emph{Complex analysis and special topics in
  harmonic analysis}.\hskip 1em plus 0.5em minus 0.4em\relax Springer-Verlag,
  New York, 1995.

\bibitem{Stepan1989Retarded}
G.~St\'{e}p\'{a}n, \emph{Retarded dynamical systems: stability and
  characteristic functions}, ser. Pitman Research Notes in Mathematics
  Series.\hskip 1em plus 0.5em minus 0.4em\relax Longman Scientific \&
  Technical, Harlow; copublished in the United States with John Wiley \& Sons,
  Inc., New York, 1989, vol. 210.

\bibitem{Wielonsky2001Rolle}
F.~Wielonsky, ``A {R}olle's theorem for real exponential polynomials in the
  complex domain,'' \emph{J. Math. Pures Appl. (9)}, vol.~80, no.~4, pp.
  389--408, 2001.

\bibitem{Amrane2018Qualitative}
S.~Amrane, F.~Bedouhene, I.~Boussaada, and S.-I. Niculescu, ``On qualitative
  properties of low-degree quasipolynomials: further remarks on the spectral
  abscissa and rightmost-roots assignment,'' \emph{Bull. Math. Soc. Sci. Math.
  Roumanie (N.S.)}, vol. 61(109), no.~4, pp. 361--381, 2018.

\bibitem{BedouheneReal}
F.~Bedouhene, I.~Boussaada, and S.-I. Niculescu, ``Real spectral values
  coexistence and their effect on the stability of time-delay systems:
  Vandermonde matrices and exponential decay,'' unpublished, available at
  \url{https://hal.archives-ouvertes.fr/hal-02476403}.

\bibitem{Boussaada2018Dominancy}
I.~Boussaada and S.-I. Niculescu, ``On the dominancy of multiple spectral
  values for time-delay systems with applications,'' \emph{IFAC-PapersOnLine},
  vol.~51, no.~14, pp. 55 -- 60, 2018, 14th IFAC Workshop on Time Delay Systems
  TDS 2018.

\bibitem{Boussaada2020Multiplicity}
I.~Boussaada, S.-I. Niculescu, A.~El-Ati, R.~P\'{e}rez-Ramos, and K.~Trabelsi,
  ``Multiplicity-induced-dominancy in parametric second-order delay
  differential equations: {A}nalysis and application in control design,''
  \emph{ESAIM Control Optim. Calc. Var.}, vol.~26, pp. Paper No. 57, 34, 2020.

\bibitem{MazantiMultiplicity}
G.~Mazanti, I.~Boussaada, and S.-I. Niculescu, ``Multiplicity-induced-dominancy
  for delay-differential equations of retarded type,'' unpublished, available
  at \url{https://hal.archives-ouvertes.fr/hal-02479909}.

\bibitem{Mazanti2020Qualitative}
------, ``On qualitative properties of single-delay linear retarded
  differential equations: Characteristic roots of maximal multiplicity are
  necessarily dominant,'' in \emph{IFAC-PapersOnLine}, 2020, in press, 21st
  IFAC World Congress.

\bibitem{Mazanti2020Spectral}
G.~Mazanti, I.~Boussaada, S.-I. Niculescu, and T.~Vyhl{\'{\i}}dal, ``Spectral
  dominance of complex roots for single-delay linear equations,'' in
  \emph{IFAC-PapersOnLine}, 2020, in press, 21st IFAC World Congress.

\bibitem{Manitius1979Finite}
A.~Z. Manitius and A.~W. Olbrot, ``Finite spectrum assignment problem for
  systems with delays,'' \emph{IEEE Trans. Automat. Control}, vol.~24, no.~4,
  pp. 541--553, 1979.

\bibitem{Avanessoff2014Hinfty}
D.~Avanessoff, A.~R. Fioravanti, C.~Bonnet, and L.~H.~V. Nguyen,
  ``{$H_\infty$}-stability analysis of (fractional) delay systems of retarded
  and neutral type with the {M}atlab toolbox {YALTA},'' in \emph{Delay
  systems}, ser. Adv. Delays Dyn.\hskip 1em plus 0.5em minus 0.4em\relax
  Springer, Cham, 2014, vol.~1, pp. 285--297.

\bibitem{Breda2009Trace}
D.~Breda, S.~Maset, and R.~Vermiglio, ``T{RACE}-{DDE}: a tool for robust
  analysis and characteristic equations for delay differential equations,'' in
  \emph{Topics in time delay systems}, ser. Lect. Notes Control Inf. Sci.\hskip
  1em plus 0.5em minus 0.4em\relax Springer, Berlin, 2009, vol. 388, pp.
  145--155.

\bibitem{Engelborghs2002Numerical}
K.~Engelborghs, T.~Luzyanina, and D.~Roose, ``Numerical bifurcation analysis of
  delay differential equations using {DDE}-{BIFTOOL},'' \emph{ACM Trans. Math.
  Software}, vol.~28, no.~1, pp. 1--21, 2002.

\bibitem{Vyhlidal2014QPmR}
T.~Vyhl\'{\i}dal and P.~Z\'{\i}tek, ``Q{P}m{R}---quasi-polynomial root-finder:
  algorithm update and examples,'' in \emph{Delay systems}, ser. Adv. Delays
  Dyn.\hskip 1em plus 0.5em minus 0.4em\relax Springer, Cham, 2014, vol.~1, pp.
  299--312.

\bibitem{Kravanja2000Computing}
P.~Kravanja and M.~Van~Barel, \emph{Computing the zeros of analytic functions},
  ser. Lecture Notes in Mathematics.\hskip 1em plus 0.5em minus 0.4em\relax
  Springer-Verlag, Berlin, 2000, vol. 1727.

\end{thebibliography}

\end{document}